# CLOSED-FORM LIKELIHOOD EXPANSIONS
# FOR MULTIVARIATE DIFFUSIONS

By Yacine Aït-Sahalia[1]

*Princeton University*

This paper provides closed-form expansions for the log-likelihood function of multivariate diffusions sampled at discrete time intervals. The coefficients of the expansion are calculated explicitly by exploiting the special structure afforded by the diffusion model. Examples of interest in financial statistics and Monte Carlo evidence are included, along with the convergence of the expansion to the true likelihood function.

**1. Introduction.** Diffusions and, more generally, continuous-time Markov processes are generally specified in economics and finance by their evolution over infinitesimal instants, that is, by writing down the stochastic differential equation followed by the state vector. However, for most estimation strategies relying on discretely sampled data, we need to be able to infer the implications of the infinitesimal time evolution of the process for the finite time intervals at which the process is actually sampled, say daily or weekly. The transition function plays a key role in that context. Unfortunately, the transition function is, in most cases, unknown.

At the same time, continuous-time models in finance, which until recently have been largely univariate, now predominantly include multiple state variables. Typical examples include asset pricing models with multiple explanatory factors, term structure models with multiple yields or factors and stochastic volatility or stochastic mean reversion models (see Sundaresan [28] for a recent survey). Motivated by this trend and the need for effective representation methods, I construct closed-form expansions for the log-transition function of a large class of multivariate diffusions. Because diffusions are Markov processes, the log-likelihood function of observations from such a process sampled at finite time intervals reduces to the sum of

Received May 2004; revised January 2007.
[1]Supported in part by NSF Grants SES-03-50772 and DMS-05-32370.

the log-transition function of successive pairs of observations. A closed form expansion for the latter therefore makes quasi-likelihood inference feasible for these models.

The paper is organized as follows. Section 2 sets out the model and assumptions. In Section 3, I introduce the concept of reducibility of a diffusion and provide a necessary and sufficient condition for the reducibility of a multivariate diffusion. In an earlier work (Aït-Sahalia [2]), I constructed explicit expansions for the transition function of univariate diffusions based on Hermite series. The natural extension of the Hermite method to the multivariate case is applicable only if the diffusion is reducible, which all univariate, but few multivariate, diffusions are. So, this paper proposes an alternative method, which determines the coefficients in closed form by requiring that the expansion satisfies the Kolmogorov equations describing the evolution of the process up to the order of the expansion itself. When a diffusion is reducible, the coefficients of the expansion are obtained as a series in the time variable, which I show in Section 4. When the diffusion is not reducible, the expansion involves a double series in the time and state variables, described in Section 5. Section 6 then studies the convergence of the likelihood expansion and the resulting maximizer to the theoretical (but incomputable) maximum likelihood estimator. Section 7 contains examples of multivariate diffusions and Monte Carlo simulation results. Proofs are in Section 8 and Section 9 concludes the paper.

**2. Setup and assumptions.** Consider the multivariate time-homogenous diffusion

$$dX_t = \mu(X_t)\,dt + \sigma(X_t)\,dW_t, \tag{1}$$

where $X_t$ and $\mu(X_t)$ are $m \times 1$ vectors, $\sigma(X_t)$ is an $m \times m$ matrix and $W_t$ is an $m \times 1$ vector of independent Brownian motions. Independence is without loss of generality since arbitrary correlation structures between the shocks to the different equations can be modeled through the inclusion of off-diagonal terms in the $\sigma$ matrix, which, furthermore, need not be symmetric. In time-inhomogeneous diffusions, the coefficients are allowed to depend on time directly, as in $\mu(X_t, t)$ and $\sigma(X_t, t)$, beyond their dependence on time via the state vector. The time-inhomogeneous case can be reduced to the time-homogenous case by treating time as an additional state variable and so it suffices to return to the model specified in (1).

The objective of this paper is to derive closed-form approximations to the log of the transition function $p_X(x|x_0, \Delta)$, that is, the conditional density of $X_{t+\Delta} = x$ given $X_t = x_0$ induced by the model (1). From an inference perspective, the primary use of this construction is to make feasible the computation of the MLE. Assume that we parametrize $(\mu, \sigma)$ as functions of a parameter vector $\theta$ and observe $X$ at dates $\{t = i\Delta \,|\, i = 0, \ldots, n\}$, where

$\Delta > 0$ is fixed. The Markovian nature of (1), which the discrete data inherit, implies that the log-likelihood function has the simple form

$$\ell_n(\theta, \Delta) \equiv \sum_{i=1}^{n} l_X(X_{i\Delta}|X_{(i-1)\Delta}, \Delta), \tag{2}$$

where $l_X \equiv \ln p_X$ and where the asymptotically irrelevant density of the initial observation, $X_0$, has been left out. In practice, the issue is that for most models of interest, the function $p_X$, hence $l_X$, is not available in closed form.

I will use the following notation. Let $\mathcal{S}_X$, a subset of $\mathbb{R}^m$, denote the domain of the diffusion $X$, assumed, for simplicity, to be of the following form.

ASSUMPTION 1. $\mathcal{S}_X$ is a product of $m$ intervals with limits $\underline{x}_i$ and $\bar{x}_i$, where possibly $\underline{x}_i = -\infty$ and/or $\bar{x}_i = +\infty$, in which case, the intervals are open at infinite limits.

I will use $^T$ to denote transposition and, for a function $\eta(x) = (\eta_1(x), \ldots, \eta_d(x))^T$, differentiable in $x$, I will write $\nabla \eta(x)$ for the Jacobian matrix of $\eta$, that is, the matrix $\nabla \eta(x) = [\partial \eta_i(x)/\partial x_j]_{i=1,\ldots,d; j=1,\ldots,m}$. For $x \in \mathbb{R}^m$, $\|x\|$ denotes the usual Euclidean norm. If $a = [a_{ij}]_{i,j=1,\ldots,m}$ is an $m \times m$ invertible matrix, then I write $a^{-1}$ for the matrix inverse, with elements $[a^{-1}]_{ij}$. $\text{Det}[a]$ and $\text{tr}[a]$ denote the determinant of $a$ and its trace, respectively. If $a = [a_i]_{i=1,\ldots,m}$ is a vector, $\text{tr}[a]$ denotes the sum of the elements of $a$. $a = \text{diag}[a_i]_{i=1,\ldots,m}$ denotes the $m \times m$ diagonal matrix with diagonal elements $a_i$. When a function $\eta(x)$ is invertible in $x$, I write $\eta^{\text{inv}}(y)$ for its inverse.

In some instances, it may be more natural to directly parametrize the infinitesimal variance–covariance matrix of the process

$$v(x) \equiv \sigma(x)\sigma^T(x) \tag{3}$$

than $\sigma(x)$ itself. Every characterization of the process, such as its transition probability, depends, in fact, on $(\mu, v)$. In particular, it can be shown that, should there exist a continuum of solutions in $\sigma$ to equation (3), then the transition probability of the process is identical for each of these $\sigma$ (see Remark 5.1.7 and Section 5.3 in Stroock and Varadhan [27]). This is also quite clear from inspection of the infinitesimal generator $A_X$ of the process, which depends on $v$ rather than $\sigma$. For functions $f(\Delta, x)$ that are suitably differentiable on its domain, $A_X$ has the action

$$A_X \cdot f = \frac{\partial f(x, \Delta)}{\partial \Delta} + \sum_{i=1}^{m} \mu_i(x) \frac{\partial f(x, \Delta)}{\partial x_i} + \frac{1}{2} \sum_{i,j=1}^{m} v_{ij}(x) \frac{\partial^2 f(x, \Delta)}{\partial x_i \partial x_j}. \tag{4}$$

The domain of $A_X$ includes at least those functions that, for each $x_0 \in \mathcal{S}_X$, are once continuously differentiable in $\Delta$ in $\mathbb{R}_+$, twice continuously differentiable in $x \in \mathcal{S}_X$ and have compact support.

As this will play a role in the likelihood expansions, define

$$D_v(x) \equiv \tfrac{1}{2} \ln(\mathrm{Det}[v(x)]). \tag{5}$$

I will assume that this matrix $v$ satisfies the following regularity condition:

ASSUMPTION 2. *The matrix $v(x)$ is positive definite for all $x$ in the interior of $\mathcal{S}_X$.*

Further assumptions are required to ensure the existence and uniqueness of a solution to (1) and to make the computation of likelihood expansions possible. I will assume the following.

ASSUMPTION 3. *$\mu(x)$ and $\sigma(x)$ are infinitely differentiable in $x$ on $\mathcal{S}_X$.*

Assumption 3 ensures the uniqueness of solutions to (1). Indeed, Assumption 3 implies in particular, that the coefficients of the stochastic differential equation are locally Lipschitz under their assumed (once) differentiability, which can be seen by applying the mean value theorem. This ensures that a solution, if it exists, will be unique (see, e.g., Theorem 5.2.5 in Karatzas and Shreve [21]). The infinite differentiability assumption in $x$ is unnecessary for that purpose, but it allows the computation of expansions of the transition density, which, as we will see, involves repeated differentiation of the coefficient functions $\mu$ and $\sigma$. There exist models of interest in finance, such as Feller's square root diffusion used in the Cox, Ingersoll and Ross model of the term structure, that fail to satisfy the Lipschitz condition since they violate the differentiability requirement of Assumption 3 at a boundary of $\mathcal{S}_X$: for instance, $\sigma(x) = \sigma_0 x^{1/2}$ is not differentiable at the left boundary 0 of $\mathcal{S}_X$. It is then possible to weaken Assumption 3 to cover such cases (see Watanabe and Yamada [30] and Yamada and Watanabe [32]).

The next assumption restricts the growth behavior of the coefficients near the boundaries of the domain.

ASSUMPTION 4. *The drift and diffusion functions satisfy linear growth conditions, that is, there exists a constant $K$ such that for all $x \in \mathcal{S}_X$ and $i, j$,*

$$|\mu_i(x)| \leq K(1 + \|x\|) \quad \text{and} \quad |\sigma_{ij}(x)| \leq K(1 + \|x\|). \tag{6}$$

*Their derivatives exhibit at most polynomial growth.*

The role of Assumption 4 is to ensure the existence of a solution to the stochastic differential equation (1) by preventing explosions of the process in finite expected time. While it can be relaxed in specific examples, it is not possible to do so in full generality. In dimension one, however, finer results are available (see the Engelbert–Schmidt criterion in Theorem 5.5.15 of Karatzas and Shreve [21]), allowing linear growth to be imposed only when the drift coefficient pulls the process toward an infinity boundary (see Proposition 1 of Aït-Sahalia [2]). In all dimensions, the linear growth condition in Assumption 4 is only an issue near the boundaries of $\mathcal{S}_X$. In the special case where $\mathcal{S}_X$ is compact, the growth condition (boundedness, in fact) follows from the continuity of the functions. The additional assumption that the derivatives of the drift and diffusion functions grow at most polynomially simplifies matters in light of the exponential tails of the transition density $p_X$.

Finally, the diffusion process $X$ is fully defined by the specification of the functions $\mu$ and $\sigma$ *and* its behavior at the boundaries of $\mathcal{S}_X$. In many examples, the specification of $\mu$ and $\sigma$ predetermines the boundary behavior of the process, but this will not be the case for models that represent limiting situations. For instance, in Cox, Ingersoll and Ross processes with affine $\mu$ and $v$, the behavior at the 0 boundary depends upon the values of the parameters. When this situation occurs for a particular model, the behavior of the likelihood expansion near such a boundary will be specified exogenously to match that of the assumed model.

**3. Reducible diffusions.** Whenever possible, I will first transform the diffusion $X$ into one that is more amenable to the derivation of an expansion for its transition density. For that purpose, I introduce the following definition.

DEFINITION 1 (*Reducibility*). The diffusion $X$ is said to be *reducible to unit diffusion* (or *reducible*, in short) if and if only if there exists a one-to-one transformation of the diffusion $X$ into a diffusion $Y$ whose diffusion matrix $\sigma_Y$ is the identity matrix. That is, there exists an invertible function $\gamma(x)$, infinitely differentiable in $X$ on $\mathcal{S}_X$, such that $Y_t \equiv \gamma(X_t)$ satisfies the stochastic differential equation

$$dY_t = \mu_Y(Y_t)\,dt + dW_t \tag{7}$$

on the domain $\mathcal{S}_Y$.

By Itô's lemma, when the diffusion is reducible, the change of variable $\gamma$ satisfies $\nabla \gamma(x) = \sigma^{-1}(x)$. Every scalar (i.e., one-dimensional) diffusion is reducible, by means of the simple transformation

$$Y_t \equiv \gamma(X_t) = \int^{X_t} \frac{du}{\sigma(u)}, \tag{8}$$

where the lower bound of integration is an arbitrary point in the interior of $\mathcal{S}_X$. The differentiability of $\gamma$ ensures that $\mu_Y$ satisfies Assumption 3. This change of variable, known as the *Lamperti transform*, played a critical role in the derivation of closed-form Hermite approximations to the transition density of univariate diffusions in Aït-Sahalia [2]. How to deal with the case where $1/\sigma(u)$ cannot be integrated in closed form is discussed after Proposition 2 below. Whenever a diffusion is reducible, an expansion can be computed for the transition density $p_X$ of $X$ by first computing it for the density $p_Y$ of the reduced process $Y$ and then transforming $Y$ back into $X$, essentially proceeding by extending the univariate method.

However, not every multivariate diffusion is reducible. Whether or not a given multivariate diffusion is reducible depends on the specification of its $\sigma$ matrix, in the following way.

PROPOSITION 1 (Necessary and sufficient condition for reducibility). *The diffusion $X$ is reducible if and only if*

$$(9) \qquad \sum_{l=1}^{m} \frac{\partial \sigma_{ik}(x)}{\partial x_l} \sigma_{lj}(x) = \sum_{l=1}^{m} \frac{\partial \sigma_{ij}(x)}{\partial x_l} \sigma_{lk}(x)$$

*for each $x$ in $\mathcal{S}_X$ and triplet $(i,j,k) = 1, \ldots, m$ such that $k > j$. If $\sigma$ is non-singular, then the condition can be expressed as*

$$(10) \qquad \frac{\partial [\sigma^{-1}]_{ij}(x)}{\partial x_k} = \frac{\partial [\sigma^{-1}]_{ik}(x)}{\partial x_j}.$$

Similar restrictions on the $\sigma$ matrix arise in different contexts; see Doss [11] who studied the question of when the solution $X$ of the SDE can be expressed as a function of the Brownian motion $W$ and the solution of an ODE and the concept of "commutative noise" in Section 10.6 of Cyganowski, Kloeden and Ombach [8] where they show that restricting the $\sigma$ matrix leads to a simplification of the Milshtein scheme for $X$.

In the bivariate case $m = 2$, condition (10) reduces to

$$\frac{\partial [\sigma^{-1}]_{11}(x)}{\partial x_2} - \frac{\partial [\sigma^{-1}]_{12}(x)}{\partial x_1} = \frac{\partial [\sigma^{-1}]_{21}(x)}{\partial x_2} - \frac{\partial [\sigma^{-1}]_{22}(x)}{\partial x_1} = 0.$$

For instance, consider diagonal systems: if $\sigma_{12} = \sigma_{21} = 0$, then the reducibility condition becomes $\partial [\sigma^{-1}]_{11}/\partial x_2 = \partial [\sigma^{-1}]_{22}/\partial x_1 = 0$. Since $[\sigma^{-1}]_{ii} = 1/\sigma_{ii}$ in the diagonal case, reducibility is equivalent to the fact that $\sigma_{ii}$ depends only on $x_i$ for each $i = 1, 2$. This is true more generally in dimension $m$. Note that this is not the case if off-diagonal elements are present. Another set of examples is provided by the class of stochastic volatility models. Consider the two models where either

$$\sigma(x) = \begin{pmatrix} \sigma_{11}(x_2) & 0 \\ 0 & \sigma_{22}(x_2) \end{pmatrix} \quad \text{or} \quad \sigma(x) = \begin{pmatrix} a(x_1) & a(x_1)b(x_2) \\ 0 & c(x_2) \end{pmatrix}.$$

In the first model, the process is not reducible in light of the previous diagonal example, as this is a diagonal system where $\sigma_{11}$ depends on $x_2$. However, in the second, the process is reducible.

**4. Closed-form expansion for the likelihood function of a reducible diffusion.** When the diffusion is reducible, I propose two approaches to construct a sequence of explicit expansions for the log-likelihood function. The first is based on computing the coefficients of a Hermite expansion for the density of the transformed process, $p_Y$. The coefficients are found in the form of a series expansion in $\Delta$, the time separating successive observations.

The second approach takes the form of the Hermite series and determines its coefficients by solving the Kolmogorov partial differential equations which characterize the transition function $p_Y$. In both cases, given a series for $p_Y$, I obtain a series for the original object of interest, $p_X$, by reversing the change of variable and the Jacobian formula. The two approaches yield the same final series.

4.1. *Multivariate Hermite expansions.* To motivate the form of the expansion that I will propose in the multivariate case, in both the reducible and irreducible cases, consider the following natural multivariate counterpart to the univariate Hermite expansion of Aït-Sahalia [2]. Let $\phi(x)$ denote the density of the $m$-dimensional multivariate Normal distribution with mean zero and identity covariance matrix. For each vector $h = (h_1, \ldots, h_m)^T \in \mathbb{N}^m$, recall that $\text{tr}[h] = h_1 + \cdots + h_m$ and let $H_h(x)$ denote the associated Hermite polynomials, which are defined by $H_h(x) = ((-1)^{\text{tr}[h]}/\phi(x)) \partial^{\text{tr}[h]} \phi(x)/dx_1^{h_1} \cdots dx_m^{h_m}$ and can be computed explicitly to an arbitrary order $\text{tr}[h]$ (see, e.g., Chapter 5 of McCullagh [23] or Withers [31]). The polynomials are orthonormal in the sense that $\int_{\mathbb{R}^m} H_h(x) H_k(x) \phi(x) \, dx = h_1! \cdots h_m!$ if $h = k$ and 0 otherwise.

The Hermite series approximation of $p_Y$ is in the form

$$(11) \quad \breve{p}_Y^{(J)}(y|y_0, \Delta) = \Delta^{-m/2} \phi\left(\frac{y - y_0}{\Delta^{1/2}}\right) \sum_{h \in \mathbb{N}^m : \text{tr}[h] \leq J} \eta_h(\Delta, y_0) H_h\left(\frac{y - y_0}{\Delta^{1/2}}\right)$$

and the Hermite coefficients $\eta_h(\Delta, y_0)$ can be computed as in the univariate case: by orthonormality of the Hermite polynomials, the coefficients $\eta_h$ are given by the conditional expectation

$$(12) \quad \eta_h(\Delta, y_0) = \frac{1}{h_1! \cdots h_m!} E[H_h(\Delta^{-1/2}(Y_{t+\Delta} - y_0)) | Y_t = y_0].$$

This expression is then amenable to computing an expansion in $\Delta$ using the generator (4). To evaluate that conditional expectation, the deterministic

Taylor expansion

$$E_{Y_1}[f(Y_\Delta, Y_0, \Delta)|Y_0 = y_0]$$
$$= \sum_{k=0}^{K} \frac{\Delta^k}{k!} A_Y^k \cdot f(y, y_0, \delta)|_{y=y_0, \delta=0} + O(\Delta^{K+1}) \quad (13)$$

can be used, where $A_Y$ is the infinitesimal generator of the process $Y$, the function $f$ is sufficiently differentiable in $(y, \delta)$ and its iterates by application of $A_Y$ up to $K$ times remain in the domain of $A_Y$, as in Aït-Sahalia [2]. The result will be a "small-time" expansion, in the same spirit as in Azencott [4] and Dacunha-Castelle and Florens-Zmirou [9], except that the expansions given here are in closed form instead of relying on moments of functionals of Brownian bridges (which are to be computed by simulation). Replacing the unknown $\eta_h$ in (11) by their expansions in $\Delta$ to order $K$ gives rises to an expansion of $\breve{p}_Y^{(J)}$ where the coefficients are gathered in increasing powers of $\Delta$, which I denote $\breve{p}_Y^{(J,K)}$. If we gather the terms in the right-hand side of (11) according to powers of $\Delta$, we can rewrite $\breve{p}_Y^{(J,K)}$ in the form of a truncated series in $\Delta$,

$$\breve{p}_Y^{(J,K)}(y|y_0, \Delta) = \Delta^{-m/2} \phi(\Delta^{-1/2}(y - y_0)) \sum_{k=0}^{K} c_Y^{(J,k)}(y|y_0) \frac{\Delta^k}{k!}. \quad (14)$$

For the log-transition density and for any given $J$, or in the univariate case where the convergence of the Hermite series is established as $J \to \infty$, the resulting expansion has the form

$$l_Y^{(K)}(y|y_0, \Delta) = -\frac{m}{2} \ln(2\pi\Delta) + \frac{C_Y^{(-1)}(y|y_0)}{\Delta} + \sum_{k=0}^{K} C_Y^{(k)}(y|y_0) \frac{\Delta^k}{k!}, \quad (15)$$

whose coefficients $C_Y^{(k)}$, $k = -1, 0, 1, 2, \ldots, K$, are combinations of the coefficients of (11) obtained by identifying the terms in the expansion in $\Delta$ of the log of (14).

The method just described is the natural extension to the multivariate setting of the Hermite approach employed in the univariate case in Aït-Sahalia [2]. Extensions of the univariate Hermite expansion results in two different univariate directions have been recently developed for time-inhomogeneous diffusions (Egorov, Li and Hu [13]) and for models driven by Lévy processes other than Brownian motion (Schaumburg [25] and Yu [33]). DiPietro [10] has extended the methodology to make it applicable in a Bayesian setting. The Hermite method requires, however, that the diffusion be reducible since the straight Hermite expansion will not in general converge if applied to $p_X$ directly instead of $p_Y$. And as discussed above, while all univariate diffusions are reducible, so that such a $Y$ exists, not all multivariate diffusions are. This necessitates an alternative method, which I now develop.

4.2. *Connection to the Kolmogorov equations.* An alternative method to obtain an explicit expansion for $l_Y$ is to take inspiration from the form of the solution given by the expansion (15) and to use the Kolmogorov equations to determine its coefficients, without any further reference to the Hermite expansion. As is often the case when a differential operator is involved, it is easier to verify that a given functional form, in this case the expansion in the form (15), is the right solution.

Consider the forward and backward Kolmogorov equations (see, e.g., Section 5.1 of Karatzas and Shreve [21])

$$(16) \quad \frac{\partial p_Y(y|y_0, \Delta)}{\partial \Delta} = -\sum_{i=1}^{m} \frac{\partial \{\mu_{Yi}(y) p_Y(y|y_0, \Delta)\}}{\partial y_i} + \frac{1}{2} \sum_{i=1}^{m} \frac{\partial^2 p_Y(y|y_0, \Delta)}{\partial y_i^2},$$

$$(17) \quad \frac{\partial p_Y(y|y_0, \Delta)}{\partial \Delta} = \sum_{i=1}^{m} \mu_{Yi}(y_0) \frac{\partial p_Y(y|y_0, \Delta)}{\partial y_{0i}} + \frac{1}{2} \sum_{i=1}^{m} \frac{\partial^2 p_Y(y|y_0, \Delta)}{\partial y_{0i}^2}.$$

The solution $p_Y$ inherits the smoothness in $(\Delta, y, y_0)$ of the coefficients $\mu_Y$ (see, e.g., Section 9.6 in Friedman [14]), so we are entitled to look for an approximate solution in the form of a smooth expansion. The fact that the Hermite expansion turns out to have exactly the right form for solving the forward and backward equations term by term is an interesting feature of these expansions.

Focusing for now on the forward equation (16), the equivalent form for the log-likelihood $l_Y$ (which is the object of interest for MLE and which will turn out to lead to a simple linear system) is

$$(18) \quad \begin{aligned} \frac{\partial l_Y(y|y_0, \Delta)}{\partial \Delta} &= -\sum_{i=1}^{m} \frac{\partial \mu_{Yi}(y)}{\partial y_i} - \sum_{i=1}^{m} \mu_{Yi}(y) \frac{\partial l_Y(y|y_0, \Delta)}{\partial y_i} \\ &\quad + \frac{1}{2} \sum_{i=1}^{m} \frac{\partial^2 l_Y(y|y_0, \Delta)}{\partial y_i^2} + \frac{1}{2} \sum_{i=1}^{m} \left( \frac{\partial l_Y(y|y_0, \Delta)}{\partial y_i} \right)^2. \end{aligned}$$

Suppose that we substitute the postulated form of the solution (15) into (18). Since

$$\begin{cases} \dfrac{\partial l_Y^{(K)}(y|y_0, \Delta)}{\partial \Delta} = -\dfrac{C_Y^{(-1)}(y|y_0)}{\Delta^2} - \dfrac{m}{2\Delta} + \sum_{k=1}^{K-1} C_Y^{(k)}(y|y_0) \dfrac{\Delta^{k-1}}{(k-1)!} \\ \dfrac{\partial l_Y^{(K)}(y|y_0, \Delta)}{\partial y_i} = \dfrac{1}{\Delta} \dfrac{\partial C_Y^{(-1)}(y|y_0)}{\partial y_i} + \sum_{k=0}^{K} \dfrac{\partial C_Y^{(-1)}(y|y_0)}{\partial y_i} \dfrac{\Delta^k}{k!} \\ \dfrac{\partial^2 l_Y^{(K)}(y|y_0, \Delta)}{\partial y_i^2} = \dfrac{1}{\Delta} \dfrac{\partial^2 C_Y^{(-1)}(y|y_0)}{\partial y_i^2} + \sum_{k=0}^{K} \dfrac{\partial^2 C_Y^{(-1)}(y|y_0)}{\partial y_i^2} \dfrac{\Delta^k}{k!}, \end{cases}$$

equating the coefficients of $\Delta^{-2}$ on both sides of (18) implies that the leading coefficient in the expansion, $C_Y^{(-1)}$, must solve the nonlinear equation

$$(19) \qquad C_Y^{(-1)}(y|y_0) = -\frac{1}{2}\left(\frac{\partial C_Y^{(-1)}(y|y_0)}{\partial y_i}\right)^T\left(\frac{\partial C_Y^{(-1)}(y|y_0)}{\partial y_i}\right).$$

Because the density must approximate a Gaussian density as $\Delta \to 0$, the appropriate solution is the one with a strict maximum at $y = y_0$, namely

$$(20) \qquad C_Y^{(-1)}(y|y_0) = -\tfrac{1}{2}\|y - y_0\|^2 = -\tfrac{1}{2}\sum_{i=1}^{m}(y_i - y_{0i})^2.$$

Considering now the coefficients of $\Delta^{-1}$ on both sides of (18), we see that

$$\sum_{i=1}^{m}\frac{\partial C_Y^{(0)}(y|y_0)}{\partial y_i}(y_i - y_{0i}) = \sum_{i=1}^{m}\mu_{Yi}(y)(y_i - y_{0i}).$$

Integrating along a line segment between $y_0$ and $y$, we obtain

$$(21) \qquad C_Y^{(0)}(y|y_0) = \sum_{i=1}^{m}(y_i - y_{0i})\int_0^1 \mu_{Yi}(y_0 + u(y - y_0))\,du,$$

with integration constants determined in the proof of the theorem below using boundary conditions and the backward equation. The higher-order coefficients are obtained using the same principle, and we have the following result.

THEOREM 1. *The coefficients of the log-density expansion $l_Y^{(K)}(y|y_0, \Delta)$ are given explicitly by* (20), (21) *and, for $k \geq 1$,*

$$(22) \qquad C_Y^{(k)}(y|y_0) = k\int_0^1 G_Y^{(k)}(y_0 + u(y - y_0)|y_0)u^{k-1}\,du.$$

*The functions $G_Y^{(k)}$ are given by*

$$G_Y^{(1)}(y|y_0) = -\sum_{i=1}^{m}\frac{\partial \mu_{Yi}(y)}{\partial y_i} - \sum_{i=1}^{m}\mu_{Yi}(y)\frac{\partial C_Y^{(0)}(y|y_0)}{\partial y_i}$$

$$(23)$$

$$+ \frac{1}{2}\sum_{i=1}^{m}\left\{\frac{\partial^2 C_Y^{(0)}(y|y_0)}{\partial y_i^2} + \left[\frac{\partial C_Y^{(0)}(y|y_0)}{\partial y_i}\right]^2\right\}$$

*and, for $k \geq 2$,*

$$G_Y^{(k)}(y|y_0) = -\sum_{i=1}^{m}\mu_{Yi}(y)\frac{\partial C_Y^{(k-1)}(y|y_0)}{\partial y_i} + \frac{1}{2}\sum_{i=1}^{m}\frac{\partial^2 C_Y^{(k-1)}(y|y_0)}{\partial y_i^2}$$

$$(24)$$

$$+ \frac{1}{2}\sum_{i=1}^{m}\sum_{h=0}^{k-1}\binom{k-1}{h}\frac{\partial C_Y^{(h)}(y|y_0)}{\partial y_i}\frac{\partial C_Y^{(k-1-h)}(y|y_0)}{\partial y_i}.$$

Theorem 1 provides the explicit form of $l_Y^{(K)}$ that solves the Kolmogorov equations to the desired order $\Delta^K$. This does not necessarily imply that $l_Y^{(K)}$ is a proper Taylor expansion of $l_Y$ at the desired order $\Delta^{K-1}$; this will be established as part of Theorem 3 below.

4.3. *Change of variable.* Given $l_Y$, the expression for $l_X$ is given by the Jacobian formula

$$
\begin{aligned}
l_X(x|x_0, \Delta) &= -\tfrac{1}{2}\ln(\text{Det}[v(x)]) + l_Y(\Delta, \gamma(x)|\gamma(x_0)) \\
&= -D_v(x) + l_Y(\Delta, \gamma(x)|\gamma(x_0)),
\end{aligned}
\tag{25}
$$

which I mimic at the level of the approximations of order $K$ in $\Delta$, thereby defining $l_X^{(K)}$ as

$$
\begin{aligned}
l_X^{(K)}(x|x_0, \Delta) &= -D_v(x) + l_Y^{(K)}(\Delta, \gamma(x)|\gamma(x_0)) \\
&= -\frac{m}{2}\ln(2\pi\Delta) - D_v(x) \\
&\quad + \frac{C_Y^{(-1)}(\gamma(x)|\gamma(x_0))}{\Delta} + \sum_{k=0}^{K} C_Y^{(k)}(\gamma(x)|\gamma(x_0))\frac{\Delta^k}{k!}
\end{aligned}
\tag{26}
$$

from $l_Y^{(K)}$ given in (15), using the coefficients $C_Y^{(k)}, k = -1, 0, \ldots, K$, given in Theorem 1. By construction, $l_X^{(K)}$ solves the Kolmogorov equations for $X$ at the same order.

**5. Closed-form expansion for the log-likelihood function of an irreducible diffusion.** In the reducible case, the two approaches (Hermite and solution of the Kolmogorov equations) coincide. When the diffusion is irreducible, however, one no longer has the option of first transforming $X$ to $Y$, computing the Hermite expansion for $Y$ and then, via the Jacobian formula, transforming it into an expansion for $X$. But, it remains possible to postulate an appropriate form of an expansion for $l_X$ and then to determine that its coefficients satisfy the Kolmogorov equations to the relevant order, as follows.

Mimicking the form of the expansion in $\Delta$ obtained in the reducible case, namely (26), leads to the postulation of the following form for an expansion of the log-likelihood

$$
\begin{aligned}
l_X^{(K)}(x|x_0, \Delta) &= -\frac{m}{2}\ln(2\pi\Delta) - D_v(x) \\
&\quad + \frac{C_X^{(-1)}(x|x_0)}{\Delta} + \sum_{k=0}^{K} C_X^{(k)}(x|x_0)\frac{\Delta^k}{k!}
\end{aligned}
\tag{27}
$$

and solving for the coefficients using the Kolmogorov equations. The expansion exists because the log-transition function inherits the smoothness of the coefficients $(\mu, v)$ (see, e.g., Section 9.6 of Friedman [14]).

When written directly for the $X$ process, as required in the irreducible case, the equations take the form

$$
\begin{aligned}
\frac{\partial l_X(x|x_0, \Delta)}{\partial \Delta} = & -\sum_{i=1}^{m} \frac{\partial \mu_i(x)}{\partial x_i} + \frac{1}{2} \sum_{i,j=1}^{m} \frac{\partial^2 v_{ij}(x)}{\partial x_i \partial x_j} \\
& - \sum_{i=1}^{m} \mu_i(x) \frac{\partial l_X(x|x_0, \Delta)}{\partial x_i} \\
& + \sum_{i,j=1}^{m} \frac{\partial v_{ij}(x)}{\partial x_i} \frac{\partial l_X(x|x_0, \Delta)}{\partial x_j} \\
& + \frac{1}{2} \sum_{i,j=1}^{m} v_{ij}(x) \frac{\partial^2 l_X(x|x_0, \Delta)}{\partial x_i \partial x_j} \\
& + \frac{1}{2} \sum_{i,j=1}^{m} \frac{\partial l_X(x|x_0, \Delta)}{\partial x_i} v_{ij}(x) \frac{\partial l_X(x|x_0, \Delta)}{\partial x_j},
\end{aligned}
\tag{28}
$$

$$
\begin{aligned}
\frac{\partial l_X(x|x_0, \Delta)}{\partial \Delta} = & \sum_{i=1}^{m} \mu_i(x_0) \frac{\partial l_X(x|x_0, \Delta)}{\partial x_{0i}} \\
& + \frac{1}{2} \sum_{i,j=1}^{m} v_{ij}(x_0) \frac{\partial^2 l_X(x|x_0, \Delta)}{\partial x_{0i} \partial x_{0j}} \\
& + \frac{1}{2} \sum_{i,j=1}^{m} \frac{\partial l_X(x|x_0, \Delta)}{\partial x_{0i}} v_{ij}(x_0) \frac{\partial l_X(x|x_0, \Delta)}{\partial x_{0j}}.
\end{aligned}
\tag{29}
$$

The solution method is as follows: as in the reducible case, substituting the postulated solution (27) into (28) provides an equation at each order in $\Delta$ which is solved for the corresponding coefficient of the expansion. While the differential equation for $l_X$ is nonlinear, it can be transformed into a linear one by exponentiation and so the expansion $l_X^{(K)}$ constructed in this way will approximate $l_X$.

Start with the equation of order $\Delta^{-2}$ which determines the leading order coefficient $C_X^{(-1)}$. While the leading coefficient $C_X^{(-1)}$ in the case of a reducible diffusion is simply $C_X^{(-1)}(x|x_0) = -\|\gamma(x) - \gamma(x_0)\|^2/2$, the situation is more involved when the diffusion $X$ is not reducible. The equation that determines the coefficient $C_X^{(-1)}$ is obtained by equating the terms of order $\Delta^{-2}$ in (28),

yielding

$$(30) \quad C_X^{(-1)}(x|x_0) = -\frac{1}{2}\left(\frac{\partial C_X^{(-1)}(x|x_0)}{\partial x}\right)^T v(x) \left(\frac{\partial C_X^{(-1)}(x|x_0)}{\partial x}\right).$$

The solution of this equation is not explicit in general, although it has a nice geometric interpretation as minus one half the square of the shortest distance from $x$ to $x_0$ in the metric induced in $\mathbb{R}^m$ by the matrix $v(x)^{-1}$ (see [29]).

5.1. *Time and state expansion.* The analysis of the coefficient $C_X^{(-1)}$ suggests that it will generally be impossible to explicitly characterize the coefficients of the expansion (27) since (30) will not in general admit an explicit solution. This is where the next step in the analysis comes into play. The idea now is to derive an explicit approximation in $(x - x_0)$ of the coefficients $C_X^{(k)}(x|x_0)$, $k = -1, 0, \ldots, K$. In other words, I localize the log-likelihood function in both $\Delta$ and $x - x_0$. The key difference between what can be done in the reducible special case of Theorem 1 and in the general case of Theorem 2 is that the coefficients of the expansion in $\Delta$ can be obtained directly by (20)–(22) with no need for an expansion in the state variable.

How this works can be seen by once again considering the coefficient $C_X^{(-1)}$. Consider a quadratic [in $(x - x_0)$] approximation of the solution to the equation (30) determining $C_X^{(-1)}$. The constant and linear terms are necessarily zero since the matrix $v(x)$ is nonsingular. Write the second-order expansion as $C_X^{(-1)}(x|x_0) = -(1/2)(x - x_0)^T V(x - x_0) + o(\|x - x_0\|^2)$. Equation (30) implies the equation $V = Vv(x_0)V$, whose solution is $V = v^{-1}(x_0)$.

As a consequence, the leading term coming from the expansion of $C_X^{(-1)}(x|x_0)$ in $x - x_0$ is $-(1/2\Delta)(x - x_0)^T v(x_0)^{-1}(x - x_0)$ so that the leading term in the expansion for the log-density will correspond to that of a Normal with mean $x_0$ and covariance matrix $\Delta v(x_0)$.

More generally, I will derive a series in $(x - x_0)$ for each coefficient $C_X^{(k)}$, at some order $j_k$ in $(x - x_0)$. That expansion is to be denoted by $C_X^{(j_k,k)}$. One remaining question is the choice of the order $j_k$ [in $(x - x_0)$] corresponding to a given order $k$ (in $\Delta$). For that purpose, note that $X_\Delta - X_0 = O_p(\Delta^{1/2})$, so

$$(31) \quad \begin{aligned} |C_X^{(k)}(X_\Delta|X_0)\Delta^k - C_X^{(j_k,k)}(X_\Delta|X_0)\Delta^k| &= O_p(\|X_\Delta - X_0\|^{j_k}\Delta^k) \\ &= O_p(\Delta^{j_k/2+k}) \end{aligned}$$

and therefore setting $j_k/2 + k = K + 1$, that is,

$$(32) \quad j_k = 2(K + 1 - k),$$

for $k = -1, 0, \ldots, K$, will provide an approximation error due to the expansion in $(x - x_0)$ of the same order $\Delta^{K+1}$ for each of the terms in the series (27).

The resulting expansion will then be of the form

$$
\begin{aligned}
\tilde{l}_X^{(K)}(x|x_0, \Delta) &= -\frac{m}{2}\ln(2\pi\Delta) - D_v(x) \\
&\quad + \frac{C_X^{(j_{-1},-1)}(x|x_0)}{\Delta} + \sum_{k=0}^{K} C_X^{(j_k, k)}(x|x_0)\frac{\Delta^k}{k!}.
\end{aligned}
\tag{33}
$$

This double expansion [in $\Delta$ and in $(x - x_0)$] can be viewed, in probability, as an expansion in $\Delta$ only once the process is inserted in the likelihood, in light of (31). In general, the function need not be analytic at $\Delta = 0$, hence this should be interpreted strictly as a series expansion.

Finally, note that the term $D_v(x)$ which arises in the reducible case from the Jacobian transformation is independent of $\Delta$ and so could be built into the $C_X^{(0)}$ coefficient. Doing so, however, would subject it to being expanded in $x - x_0$, which is unnecessary since $D_v(x)$ is known. If $D_v(x)$ were being expanded along with $C_X^{(j_0, 0)}$, we would lose the property that $\tilde{l}_X^{(K)}$ also solves the backward equation (29) to the corresponding order in $\Delta$.

5.2. *Determination of the coefficients in the irreducible case.* What remains to be done is to explicitly compute the expansion $C_X^{(j_k, k)}$ in $x - x_0$ of each coefficient $C_X^{(k)}$. Let $i \equiv (i_1, i_2, \ldots, i_m)$ denote a vector of integers and define $I_k = \{i \equiv (i_1, i_2, \ldots, i_m) \in \mathbb{N}^m : 0 \leq \text{tr}[i] \leq j_k\}$ so that the form of $C_X^{(j_k, k)}$ is

$$
C_X^{(j_k, k)}(x|x_0) = \sum_{i \in I_k} \beta_i^{(k)}(x_0)(x_1 - x_{01})^{i_1}(x_2 - x_{02})^{i_2} \cdots (x_m - x_{0m})^{i_m}.
\tag{34}
$$

The coefficients are determined one by one, starting with the leading term $C_X^{(j_{-1}, -1)}$. Given $C_X^{(j_{-1}, -1)}$, the next term $C_X^{(j_0, 0)}$ is calculated explicitly, and so on. Based on (32), the highest-order term $(k = -1)$ is expanded to a higher degree of precision $j_{-1}$ than the successive terms. This is quite natural, given that $C_X^{(j_{-1}, -1)}$ is an input to the differential equation determining $C_X^{(j_0, 0)}$, and so on. In order to state the main result pertaining to the closed-form solutions $C_X^{(j_k, k)}$, I define the following functions of the coefficients and their derivatives:

$$
\begin{aligned}
G_X^{(0)}(x|x_0) &= \frac{m}{2} - \sum_{i=1}^{m} \mu_i(x)\frac{\partial C_X^{(-1)}(x|x_0)}{\partial x_i} \\
&\quad + \sum_{i,j=1}^{m} \frac{\partial v_{ij}(x)}{\partial x_i}\frac{\partial C_X^{(-1)}(x|x_0)}{\partial x_j}
\end{aligned}
$$

(35)
$$+ \frac{1}{2} \sum_{i,j=1}^{m} v_{ij}(x) \frac{\partial^2 C_X^{(-1)}(x|x_0)}{\partial x_i \partial x_j}$$
$$- \sum_{i,j=1}^{m} v_{ij}(x) \frac{\partial C_X^{(-1)}(x|x_0)}{\partial x_i} \frac{\partial D_v(x)}{\partial x_j},$$

$$G_X^{(1)}(x|x_0) = -\sum_{i=1}^{m} \frac{\partial \mu_i(x)}{\partial x_i} + \frac{1}{2} \sum_{i,j=1}^{m} \frac{\partial^2 v_{ij}(x)}{\partial x_i \partial x_j}$$
$$- \sum_{i=1}^{m} \mu_i(x) \left( \frac{\partial C_X^{(0)}(x|x_0)}{\partial x_i} - \frac{\partial D_v(x)}{\partial x_i} \right)$$
$$+ \sum_{i,j=1}^{m} \frac{\partial v_{ij}(x)}{\partial x_i} \partial x_i \left( \frac{\partial C_X^{(0)}(x|x_0)}{\partial x_j} - \frac{\partial D_v(x)}{\partial x_j} \partial x_j \right)$$

(36)
$$+ \frac{1}{2} \sum_{i,j=1}^{m} v_{ij}(x) \left\{ \frac{\partial^2 C_X^{(0)}(x|x_0)}{\partial x_i \partial x_j} - \frac{\partial^2 D_v(x)}{\partial x_i \partial x_j} \right.$$
$$+ \left( \frac{\partial C_X^{(0)}(x|x_0)}{\partial x_i} - \frac{\partial D_v(x)}{\partial x_i} \right)$$
$$\left. \times \left( \frac{\partial C_X^{(0)}(x|x_0)}{\partial x_j} - \frac{\partial D_v(x)}{\partial x_j} \right) \right\}$$

and, for $k \geq 2$,

$$G_X^{(k)}(x|x_0) = -\sum_{i=1}^{m} \mu_i(x) \frac{\partial C_X^{(k-1)}(x|x_0)}{\partial x_i} + \sum_{i,j=1}^{m} \frac{\partial v_{ij}(x)}{\partial x_i} \frac{\partial C_X^{(k-1)}(x|x_0)}{\partial x_j}$$

(37)
$$+ \frac{1}{2} \sum_{i,j=1}^{m} v_{ij}(x) \frac{\partial^2 C_X^{(k-1)}(x|x_0)}{\partial x_i \partial x_j}$$
$$+ \frac{1}{2} \sum_{i,j=1}^{m} v_{ij}(x) \left\{ \left( \frac{\partial C_X^{(0)}(x|x_0)}{\partial x_i} - 2\frac{\partial D_v(x)}{\partial x_i} \right) \frac{\partial C_X^{(k-1)}(x|x_0)}{\partial x_j} \right.$$
$$+ \sum_{h=0}^{k-2} \binom{k-1}{h}$$
$$\left. \times \frac{\partial C_X^{(h)}(x|x_0)}{\partial x_i} \frac{\partial C_X^{(k-1-h)}(x|x_0)}{\partial x_j} \right\}.$$

Note that the computation of each function $G_X^{(k)}$ requires only the ability to differentiate the previously determined coefficients $C_X^{(-1)}, \ldots, C_X^{(k-1)}$. The same applies to their expansions. The following theorem can now describe how the coefficients $C_X^{(j_k,k)}$, that is, the coefficients $\beta_i^{(k)}$, $i \in I_k$, are determined.

THEOREM 2. *For each $k = -1, 0, \ldots, K$, the coefficient $C_X^{(k)}(x|x_0)$ in (27) solves the equation*

$$f_X^{(k-1)}(x|x_0) = 0, \tag{38}$$

*where*

$$f_X^{(-2)}(x|x_0) = -2C_X^{(-1)}(x|x_0) - \sum_{i,j=1}^{m} v_{ij}(x) \frac{\partial C_X^{(-1)}(x|x_0)}{\partial x_i} \frac{\partial C_X^{(-1)}(x|x_0)}{\partial x_j}, \tag{39}$$

$$f_X^{(-1)}(x|x_0) = -\sum_{i,j=1}^{m} v_{ij}(x) \frac{\partial C_X^{(-1)}(x|x_0)}{\partial x_i} \frac{\partial C_X^{(0)}(x|x_0)}{\partial x_j} - G_X^{(0)}(x|x_0) \tag{40}$$

*and, for $k \geq 1$,*

$$\begin{aligned}f_X^{(k-1)}(x|x_0) = {} & C_X^{(k)}(x|x_0) \\ & - \frac{1}{k} \sum_{i,j=1}^{m} v_{ij}(x) \frac{\partial C_X^{(-1)}(x|x_0)}{\partial x_i} \frac{\partial C_X^{(k)}(x|x_0)}{\partial x_j} - G_X^{(k)}(x|x_0),\end{aligned} \tag{41}$$

*where the functions $G_X^{(k)}$, $k = 0, 1, \ldots, K$, are given above. The coefficients $\beta_i^{(k)}$ solve a system of linear equations, whose solution is explicit.*

$G_X^{(k)}$ involves only the coefficients $C_X^{(h)}$ for $h = -1, \ldots, k-1$, so this system of equations can be utilized to solve recursively for each coefficient, meaning that the equation $f_X^{(-2)} = 0$ determines $C_X^{(-1)}$; given $C_X^{(-1)}$, $G_X^{(0)}$ becomes known and the equation $f_X^{(-1)} = 0$ determines $C_X^{(0)}$; given $C_X^{(-1)}$ and $C_X^{(0)}$, $G_X^{(1)}$ becomes known and the equation $f_X^{(0)} = 0$ then determines $C_X^{(1)}$, and so on.

Each of these equations can be solved explicitly in the form of the expansion $C_X^{(j_k,k)}$ of the coefficient $C_X^{(k)}$, at order $j_k$ in $(x - x_0)$. The coefficients $\beta_i^{(k)}(x_0)$, $i \in I_k$, of $C_X^{(j_k,k)}$ are determined by setting the expansion $f_X^{(j_k,k-1)}$ of $f_X^{(k-1)}$ to zero. The key feature that makes this problem solvable in closed form is that the coefficients solve a succession of systems of linear equations: first determine $\beta_i^{(k)}$ for $\text{tr}[i] = 0$, then $\beta_i^{(k)}$ for $\text{tr}[i] = 1$ and so on, all the way to $\mathrm{tr}[i] = j_k$. Note, in particular, for $k = -1$, $\beta_i^{(-1)} = 0$ for $\mathrm{tr}[i] = 0, 1$ (i.e., the polynomial has no constant or linear terms) and the terms corresponding to $\mathrm{tr}[i] = 2$ (with, of course, $j_{-1} \geq 2$) are

$$\sum_{i \in I_{-1} : \mathrm{tr}[i]=2} \beta_i^{(-1)}(x_0)(x_1 - x_{01})^{i_1} \cdots (x_m - x_{0m})^{i_m}$$

$$= -\tfrac{1}{2}(x - x_0)^T v^{-1}(x_0)(x - x_0),$$

which are the anticipated terms, given the Gaussian limiting behavior of the transition density when $\Delta$ is small. Thus, with $j_{-1} \geq 3$, we only need to determine the terms $\beta_i^{(-1)}$ corresponding to $\mathrm{tr}[i] = 3, \ldots, j_{-1}$. Note that the solution $\beta_i^{(-1)}$ depends only on the specification of the $v$ matrix (the drift functions are irrelevant). For $k = 0$, $\beta_i^{(0)} = 0$ for $\mathrm{tr}[i] = 0$, so the polynomial has no constant term. For $k \geq 1$, the polynomials have a constant term (for $k \geq 1$, $\beta_i^{(k)} \neq 0$ for $\mathrm{tr}[i] = 0$ in general).

To obtain an expansion for the density $p_X$ instead of for the log-density $l_X$, one can either take the exponential of $\tilde{l}_X^{(K)}$ or, alternatively, expand the exponential in $\Delta$ to obtain the coefficients $c_X$ for the expansion of $p_X$ from the coefficients $C_X$ for the expansion of $l_X$. In general, like a Hermite expansion, neither will integrate to one without division by the integral over $\mathcal{S}_X$ of the density expansion. Positivity is guaranteed, however, if one simply exponentiates the log-transition function.

5.3. *Applying the irreducible method to a reducible diffusion.* Theorem 2 is more general than Theorem 1, in that it does not require that the diffusion be reducible. As discussed above, in exchange for that generality, the coefficients are available in closed form only in the form of a series expansion in $x$ about $x_0$. The following proposition describes the relationship between these two methods when Theorem 2 is applied to a diffusion that is, in fact, reducible.

PROPOSITION 2. *Suppose that the diffusion $X$ is reducible and let $l_X^{(K)}$ denote its log-likelihood expansion calculated by applying Theorem 1. Suppose, now, that we also calculate its log-likelihood expansion, $\tilde{l}_X^{(K)}$, without first transforming $X$ into the unit diffusion $Y$, that is, by applying Theorem 2 to $X$ directly. Then, each coefficient $C_X^{(j_k,k)}(x|x_0)$ from $\tilde{l}_X^{(K)}$ is an expansion in $(x - x_0)$ at order $j_k$ of the coefficient $C_X^{(k)}(x|x_0) = C_Y^{(k)}(\gamma(x)|\gamma(x_0))$ from $l_X^{(K)}$.*

In other words, applying the irreducible method to a diffusion that is, in fact, reducible involves replacing the exact expression for $C_X^{(k)}(x|x_0)$ by its series in $(x - x_0)$. Of course, there is no reason to do so when the diffusion is reducible and the transformation $\gamma$ from $X$ to $Y$, defined in Definition 1, is explicit.

However, Proposition 2 is relevant in the case where the diffusion is reducible, but the transformation $\gamma$ is not available in closed form. This can occur even in dimension $m = 1$, where every diffusion is theoretically reducible. For instance, consider the specification of the diffusion function in the general interest rate model proposed in [1], namely $\sigma^2(x) = \theta_{-1}x^{-1} + \theta_0 + \theta_1 x + \theta_2 x^{\theta_3}$, where the $\theta$'s denote parameters. In that case, $\gamma(x)$, given in (8), involves integrating $1/\sigma$ and the result is not explicit. Fortunately, one can use the irreducible method in that case and the result of applying that method is given by Proposition 2. An alternative is to use the method that has been proposed by [5].

Finally, the double series in $\Delta$ and $(x - x_0)$ produced by the irreducible method matches, when applied to a reducible diffusion, the expansion produced by the Hermite series since the coefficients of the latter [a polynomial in $(x - x_0)$, by construction] are obtained as a series in $\Delta$ by computing their conditional expectations, as described in (13). But the infinitesimal generator of the process in (13) is by definition, such that the resulting coefficients solve, at each order in $\Delta$, the Kolmogorov equations. Hence, the two series match.

**6. Convergence to the true log-likelihood function and the resulting approximate MLE.** Theorems 1 and 2 give the expressions of the coefficients of the expansion in the reducible and irreducible cases, respectively. I now turn to the convergence of the resulting expansion to the object of interest, showing that the series constructed above is a Taylor expansion of the true, but unknown, log-likelihood function, and considering its application to likelihood inference.

Suppose that $(\mu, \sigma)$ are parametrized using a parameter vector $\theta$ and that $(\mu, \sigma)$ and their derivatives at all orders are three times continuously differentiable in $\theta$. The differentiability of the coefficients extends to $l_X$, in light of the previously cited results on the solutions of second-order parabolic partial differential equations (Section 9.6 of Friedman [14]), and to the expansion by construction, given that it consists of sums and products of $(\mu, \sigma)$ and their derivatives. Let the parameter space $\Theta$ be a compact subset of $\mathbb{R}^r$. Let $\theta_0$ denote the true value of the parameter. Assume, for simplicity, that for fixed $n$ and $\Delta$, $\theta \longmapsto \ell_n(\theta, \Delta)$ defined in (2) has a unique maximizer $\hat{\theta}_{n,\Delta} \in \Theta$. $\hat{\theta}_{n,\Delta}$ is the exact (but incomputable) MLE for $\theta$. Consider, now, the approximate MLE $\hat{\theta}^{(K)}_{n,\Delta}$ obtained by maximizing $\ell^{(K)}_n(\theta, \Delta)$ (resp. $\tilde{\ell}^{(K)}_n$), itself defined analogously to (2), but with the expansion $l^{(K)}_X$ (resp. $\tilde{l}^{(K)}_X$) in the reducible (resp. irreducible) case instead of the true log-transition function $l_X$. We have the following result.

THEOREM 3. *For any $n$,*

$$\sup_{\theta \in \Theta} |\tilde{\ell}_n^{(K)}(\theta, \Delta) - \ell_n(\theta, \Delta)| \to 0 \tag{42}$$

*in $P_{\theta_0}$-probability as $\Delta \to 0$. In the reducible case, the same holds for $\ell_n^{(K)}$. The approximate MLE sequence $\hat{\theta}_{n,\Delta}^{(K)}$ exists almost surely and satisfies $\hat{\theta}_{n,\Delta}^{(K)} - \hat{\theta}_{n,\Delta} \to 0$ in $P_{\theta_0}$-probability as $\Delta \to 0$.*

*Furthermore, suppose that as $\to \infty$, we have $\hat{\theta}_{n,\Delta} \to \theta_0$ in $P_{\theta_0}$-probability and that there exists a sequence of nonsingular $r \times r$ matrices $S_{n,\Delta}$ such that*

$$S_{n,\Delta}^{-1}(\hat{\theta}_{n,\Delta} - \theta_0) = O_p(1). \tag{43}$$

*There then exists a sequence $\Delta_n \to 0$ such that*

$$S_{n,\Delta_n}^{-1}(\hat{\theta}_{n,\Delta_n}^{(K)} - \hat{\theta}_{n,\Delta_n}) = o_p(1). \tag{44}$$

Intuitively, the reason that the log-approximation error (42) is small in probability is as follows. For small $\Delta$, in a small neighborhood about $x_0$, the approximation error is small by construction because $l_X^{(K)}$ (resp. $\tilde{l}_X^{(K)}$) is a Taylor expansion of $l_X$ about $\Delta = 0$ (and about $x = x_0$, resp.). Away from $x_0$, the approximation error may not be small, unless $l_X$ is analytic, but this does not matter much because such an error is at most polynomial, while the probability of reaching this region in time $\Delta$ is exponentially small.

Note, also, that it follows from (43)–(44) that $\hat{\theta}_{n,\Delta}^{(K)}$ and $\hat{\theta}_{n,\Delta}$ share the same asymptotic distribution as $\to \infty$. For instance, (43) is verified, in particular, if the process $X$ is stationary with positive definite Fisher information matrix $F_\Delta$ for a pair of successive observations, in which case $S_{n,\Delta}$ can be taken to be $n^{-1/2} F_\Delta^{1/2}$ (see Billingsley [6] for the required regularity conditions).

**7. Examples.** In this section, I apply the above results to a leading multivariate diffusion example. The last example shows that the method of this paper also applies to time-inhomogeneous models.

7.1. *The Bivariate Ornstein–Uhlenbeck model.* Consider the model

$$dX_t = \beta(\alpha - X_t)\,dt + \sigma\,dW_t, \tag{45}$$

where $\alpha = [\alpha_i]_{i=1,2}$, $\beta = [\beta_i]_{i=1,2}$ and $\sigma = [\sigma_{i,j}]_{i,j=1,2}$ and assume that $\beta$ and $\sigma$ have full rank. This is the most basic model capturing mean reversion in the state variables. Consider, now, the matrix equation $\beta\lambda + \lambda\beta^T = \sigma\sigma^T$, whose solution in the bivariate case is the $2 \times 2$ symmetric matrix $\lambda$ given by

$$\lambda = \frac{1}{2\operatorname{tr}[\beta]\operatorname{Det}[\beta]}(\operatorname{Det}[\beta]\sigma\sigma^T + (\beta - \operatorname{tr}[\beta])\sigma\sigma^T(\beta - \operatorname{tr}[\beta])^T). \tag{46}$$

When the process is stationary, that is, when the eigenvalues of the matrix $\beta$ have positive real parts, $\lambda$ is the stationary variance–covariance matrix of the process. That is, the stationary density of $X$ is the bivariate Normal density with mean $\alpha$ and variance–covariance $\lambda$.

The transition density of $X$ is the bivariate Normal density

$$
\begin{aligned}
(47)\quad p_X(x|x_0, \Delta) &= (2\pi)^{-1} \operatorname{Det}[\Omega(\Delta)]^{-1/2} \\
&\quad \times \exp(-(x - m(\Delta, x_0))^T \Omega^{-1}(\Delta)(x - m(\Delta, x_0))),
\end{aligned}
$$

where $m(\Delta, x_0) = \alpha + \exp(-\beta\Delta)(x_0 - \alpha)$ and $\Omega(\Delta) = \lambda - \exp(-\beta\Delta)\lambda\exp(-\beta^T\Delta)$, with exp denoting the matrix exponential.

Identification of the continuous-time parameters from the discrete data for this particular model is discussed in Philips [24], Hansen and Sargent [15] and Kessler and Rahbek [22]. If we wish to identify the parameters in $\theta$ from discrete data sampled at the given time interval $\Delta$, then we must restrict the set of admissible parameter values $\Theta$. For instance, we may restrict $\Theta$ in such a way that the mapping $\beta \mapsto \exp(-\beta\Delta)$ is invertible, for instance, by restricting the admissible parameter matrices $\beta$ to have real eigenvalues. This will be the case, for example, if we restrict attention to matrices $\beta$ which are triangular (and, of course, have real elements). For the rest of this discussion, I will assume that $\Theta$ has been restricted in such a way.

By applying Proposition 1, we see that the process $X$ is reducible and that $\gamma(x) = \sigma^{-1}x$, so

$$
(48)\quad dY_t = (\sigma^{-1}\beta\alpha - \sigma^{-1}\beta\sigma Y_t)\,dt + dW_t \equiv \kappa(\eta - Y_t)\,dt + dW_t,
$$

where $\eta = \sigma^{-1}\alpha = [\eta_i]_{i=1,2}$ and $\kappa = \sigma^{-1}\beta\sigma = [\kappa_{i,j}]_{i,j=1,2}$. One can therefore apply Theorem 1 to obtain the coefficients of the expansion:

$$
\begin{aligned}
C_Y^{(-1)}(y|y_0) &= -\tfrac{1}{2}(y_1 - y_{01})^2 - \tfrac{1}{2}(y_2 - y_{02})^2, \\
C_Y^{(0)}(y|y_0) &= -\tfrac{1}{2}(y_1 - y_{01})((y_1 + y_{01} - 2\gamma_1)\kappa_{11} + (y_2 + y_{02} - 2\gamma_2)\kappa_{12}) \\
&\quad - \tfrac{1}{2}(y_2 - y_{02})((y_1 + y_{01} - 2\gamma_1)\kappa_{21} + (y_2 + y_{02} - 2\gamma_2)\kappa_{22}), \\
C_Y^{(1)}(y|y_0) &= \tfrac{1}{2}(\kappa_{11} - ((y_{01} - \eta_1)\kappa_{11} + (y_{02} - \eta_2)\kappa_{12})^2) \\
&\quad + \tfrac{1}{2}(\kappa_{22} - ((y_{01} - \eta_1)\kappa_{21} + (y_{02} - \eta_2)\kappa_{22})^2) \\
&\quad - \tfrac{1}{2}(y_1 - y_{01})((y_{01} - \eta_1)(\kappa_{11}^2 + \kappa_{21}^2) \\
&\quad\quad + (y_{02} - \eta_2)(\kappa_{11}\kappa_{12} + \kappa_{21}\kappa_{22})) \\
&\quad + \tfrac{1}{24}(y_1 - y_{01})^2(-4\kappa_{11}^2 + \kappa_{12}^2 - 2\kappa_{12}\kappa_{21} - 3\kappa_{21}^2) \\
&\quad - \tfrac{1}{2}(y_2 - y_{02})((y_{01} - \eta_1)(\kappa_{11}\kappa_{12} + \kappa_{21}\kappa_{22}) \\
&\quad\quad + (y_{02} - \eta_2)(\kappa_{12}^2 + \kappa_{22}^2))
\end{aligned}
$$

$$+ \tfrac{1}{24}(y_2 - y_{02})^2(-4\kappa_{22}^2 + \kappa_{21}^2 - 2\kappa_{12}\kappa_{21} - 3\kappa_{12}^2)$$
$$- \tfrac{1}{3}(y_1 - y_{01})(y_2 - y_{02})(\kappa_{11}\kappa_{12} + \kappa_{21}\kappa_{22}),$$
$$C_Y^{(2)}(y|y_0) = -\tfrac{1}{12}(2\kappa_{11}^2 + 2\kappa_{22}^2 + (\kappa_{12} + \kappa_{21})^2)$$
$$+ \tfrac{1}{6}(y_1 - y_{01})(\kappa_{12} - \kappa_{21})((y_{01} - \eta_1)(\kappa_{11}\kappa_{12} + \kappa_{21}\kappa_{22})$$
$$+ (y_{02} - \eta_2)(\kappa_{12}^2 + \kappa_{22}^2))$$
$$+ \tfrac{1}{12}(y_1 - y_{01})^2(\kappa_{12} - \kappa_{21})(\kappa_{11}\kappa_{12} + \kappa_{21}\kappa_{22})$$
$$+ \tfrac{1}{12}(y_2 - y_{02})^2(\kappa_{21} - \kappa_{12})(\kappa_{11}\kappa_{12} + \kappa_{21}\kappa_{22})$$
$$+ \tfrac{1}{6}(y_2 - y_{02})(\kappa_{21} - \kappa_{12})((y_{01} - \eta_1)(\kappa_{11}^2 + \kappa_{21}^2)$$
$$+ (y_{02} - \eta_2)(\kappa_{11}\kappa_{12} + \kappa_{21}\kappa_{22}))$$
$$+ \tfrac{1}{12}(y_1 - y_{01})(y_2 - y_{02})(\kappa_{12} - \kappa_{21})(\kappa_{22}^2 + \kappa_{12}^2 - \kappa_{11}^2 + \kappa_{21}^2).$$

Because this is essentially the only multivariate model with a known closed-form density (other than multivariate models which reduce to the superposition of univariate processes), the Ornstein–Uhlenbeck process can serve as a useful benchmark for examining the accuracy of the expansions and the resulting MLE. Table 1 reports the results of 1,000 Monte Carlo simulations comparing the distribution of the maximum-likelihood estimator $\hat{\theta}^{(\text{MLE})}$ based on the exact transition density for this model, around the true value of the parameters $\theta^{(\text{TRUE})}$, to the distribution of the difference between the MLE $\hat{\theta}^{(\text{MLE})}$ and the approximate MLE $\hat{\theta}^{(2)}$ based on the expansion with $K = 2$ terms shown above. To ensure full identification, the off-diagonal term $\kappa_{21}$ is constrained to be zero. As discussed above, this guarantees that the eigenvalues of the mean reversion matrix are both real and avoids the aliasing problem altogether. The constraints $\kappa_{11} > 0$ and $\kappa_{22} > 0$ make the process stationary, so standard asymptotics give the asymptotic distribution of $\hat{\theta}^{(\mathrm{MLE})}$ (the inverse of Fisher's information is computed as $E[-\partial^2 l_X/\partial\theta\,\partial\theta^T]^{-1}$).

TABLE 1
*Monte Carlo simulations for the bivariate Ornstein–Uhlenbeck model*

| Parameter | $\theta^{(\text{TRUE})}$ | Asymptotic $\hat{\theta}^{(\text{MLE})} - \theta^{(\text{TRUE})}$ | | Small sample $\hat{\theta}^{(\text{MLE})} - \theta^{(\text{TRUE})}$ | | Small sample $\hat{\theta}^{(\text{MLE})} - \hat{\theta}^{(2)}$ | |
|---|---|---|---|---|---|---|---|
| | | Mean | Stdev | Mean | Stdev | Mean | Stdev |
| $\eta_1$ | 0 | 0 | 0.065 | −0.0013 | 0.066 | −0.0000005 | 0.000014 |
| $\eta_2$ | 0 | 0 | 0.032 | −0.001 | 0.033 | −0.0000003 | 0.000011 |
| $\kappa_{11}$ | 5 | 0 | 1.03 | 0.49 | 1.11 | 0.012 | 0.008 |
| $\kappa_{12}$ | 1 | 0 | 1.51 | 0.12 | 1.64 | 0.010 | 0.016 |
| $\kappa_{22}$ | 10 | 0 | 1.44 | 0.33 | 1.46 | 0.068 | 0.029 |

Each of the 1,000 samples is a series of $n = 500$ weekly observations ($\Delta = 1/52$), generated using the exact discretization of the process. The results in the table show that the difference $\hat{\theta}^{(\mathrm{MLE})} - \hat{\theta}^{(2)}$ is an order of magnitude smaller than the (inescapable) sampling error $\hat{\theta}^{(\mathrm{MLE})} - \theta^{(\mathrm{TRUE})}$. Hence, for the purpose of estimating $\theta^{(\mathrm{TRUE})}$, $\hat{\theta}^{(2)}$ can be taken as a useful substitute for the (generally incomputable) $\hat{\theta}^{(\mathrm{MLE})}$. In other words, at least for this model, $K = 2$ provides sufficient accuracy for the types of situations and values of the sampling interval $\Delta$ one typically encounters in finance.

7.2. *Comparing the accuracy of the reducible and irreducible methods.* Using nonlinear transformations of the Ornstein–Uhlenbeck process, we can assess the empirical performance of the general method for irreducible diffusions. Let $Y$ denote the process given in (48) and define $X_t \equiv \exp(Y_t) = \gamma^{\mathrm{inv}}(Y_t)$. From Itô's lemma, the dynamics of $X_t$ are given by

$$
(49) \quad dX_t = \begin{pmatrix} X_{1t}(\frac{1}{2} + \kappa_{11}(\eta_1 - \ln(X_{1t})) + \kappa_{12}(\eta_2 - \ln(X_{2t}))) \\ X_{2t}(\frac{1}{2} + \kappa_{21}(\eta_1 - \ln(X_{1t})) + \kappa_{22}(\eta_2 - \ln(X_{2t}))) \end{pmatrix} dt + \begin{pmatrix} X_{1t} & 0 \\ 0 & X_{2t} \end{pmatrix} dW_t.
$$

By construction, this process has a known log-transition function given by $l_X(x|x_0, \Delta) = -\ln(xx_0) + l_Y(\Delta, \ln(x)|\ln(x_0))$ and it is reducible by transforming $X_t$ back to $Y_t = \ln(X_t) = \gamma(X_t)$, with $D_v(x) = \ln(X_{1t}X_{2t})$ for that transformation.

But, in order to assess the accuracy of the irreducible method, we can directly calculate the irreducible expansion (based on Theorem 2) for the model (49). We can then compare it to the closed-form solution, but also to the reducible expansion obtained using the order 2 expansion given in the previous section for $l_Y$ and then the Jacobian formula, $l_X(x|x_0, \Delta) = -\ln(X_{1t}X_{2t}) + l_Y(\Delta, \ln(x)|\ln(x_0))$. Based on Proposition 2, we know that in this situation, the irreducible expansion involves Taylor expanding the coefficients of the reducible expansion in $x$ about $x_0$. Monte Carlo simulations with the same design as in the previous section help document the effect of that further Taylor expansion on the accuracy of the resulting MLE. The results are given in Table 2 and they show that the difference $\hat{\theta}^{(\mathrm{MLE})} - \hat{\theta}^{(2,\mathrm{irreducible})}$, although larger than $\hat{\theta}^{(\mathrm{MLE})} - \hat{\theta}^{(2,\mathrm{reducible})}$, remains smaller than the difference $\hat{\theta}^{(\mathrm{MLE})} - \theta^{(\mathrm{TRUE})}$ due to the sampling noise. In other words, replacing $\hat{\theta}^{(\mathrm{MLE})}$ by $\hat{\theta}^{(2,\mathrm{irreducible})}$ has an effect which is not statistically discernible from the sampling variation of $\hat{\theta}^{(\mathrm{MLE})}$ around $\theta^{(\mathrm{TRUE})}$. And, of course, $\hat{\theta}^{(\mathrm{MLE})}$ is generally incomputable, whereas $\hat{\theta}^{(2,\mathrm{irreducible})}$ is computable.

TABLE 2
*Monte Carlo simulations for the exponential transformation of the Ornstein–Uhlenbeck model: comparing the reducible and irreducible methods*

| Parameter | $\theta^{(\text{TRUE})}$ | Small sample $\hat{\theta}^{(\text{MLE})} - \theta^{(\text{TRUE})}$ | | Small sample $\hat{\theta}^{(\text{MLE})} - \hat{\theta}^{(2,\text{reducible})}$ | | Small sample $\hat{\theta}^{(\text{MLE})} - \hat{\theta}^{(2,\text{irreducible})}$ | |
|---|---|---|---|---|---|---|---|
| | | **Mean** | **Stdev** | **Mean** | **Stdev** | **Mean** | **Stdev** |
| $\eta_1$ | 0 | $-0.0019$ | 0.065 | $-0.000003$ | 0.00009 | $-0.0004$ | 0.0002 |
| $\eta_2$ | 0 | $-0.003$ | 0.034 | $-0.0000001$ | 0.00001 | 0.0003 | 0.0002 |
| $\kappa_{11}$ | 5 | 0.48 | 1.06 | 0.012 | 0.007 | 0.08 | 0.04 |
| $\kappa_{12}$ | 1 | $-0.05$ | 1.57 | 0.0087 | 0.016 | 0.025 | 0.03 |
| $\kappa_{22}$ | 10 | 0.39 | 1.50 | 0.069 | 0.030 | 0.21 | 0.08 |

7.3. *Time-inhomogeneous models.* Time-inhomogeneous models are of particular interest for the term structure of interest rates. A large swathe of the term structure literature has proposed models designed to fit exactly the current bond prices, as well as other market data, such as bond volatilities or the implied volatilities of interest rate caps, for instance. Calibrating such a model to time-varying market data gives rise to time-varying drift and diffusion coefficients. Typical examples of this approach include the models of Ho and Lee [17], Black, Derman and Toy [7] and Hull and White [18], where the short-term interest rate (or its log) follows the dynamics $dX_{1t} = (\alpha(t) - \beta(t)X_{1t})\,dt + \kappa(t)\,dW_{1t}$. Markovian specializations of the Heath, Jarrow and Morton [16] model will also be, in general, time-inhomogeneous.

The univariate results of Aït-Sahalia [2] have been extended to cover such models by Egorov, Li and Xu [13]. With expansions now available for time-homogenous diffusions of arbitrary specifications and dimensions, a time-inhomogeneous diffusion of dimension $m$ can be simply reduced to a time-homogenous diffusion of dimension $m+1$. Indeed, consider the state vector $X_t = (X_{1t}, \ldots, X_{mt})$. Now, define time as the additional state variable $X_{m+1,t} = t$, whose dynamics are $dX_{m+1,t} = dt$, and consider the extended state vector as $\tilde{X}_t = (X_{1t}, \ldots, X_{mt}, X_{m+1,t})$. This is an $(m+1)$-dimensional, time-homogenous, diffusion.

## 8. Proofs.

8.1. *Proof of Proposition* 1. Suppose that a transformation exists and define $Y_t \equiv \gamma(X_t)$, where $\gamma(x) = (\gamma_1(x), \ldots, \gamma_m(x))^T$. By Itô's lemma, the diffusion matrix of $Y$ is $\sigma_Y(Y_t) = \nabla \gamma(X_t)\,\sigma(X_t)$. For $\sigma_Y$ to be *Id*, we must therefore have that $\nabla \gamma(X_t) = \sigma^{-1}(x)$ (recall that $\sigma$ is assumed to be non-singular). Thus,

$$[\sigma^{-1}]_{ij}(x) = \frac{\partial \gamma_i(x)}{\partial x_j} \tag{50}$$

and hence

$$\frac{\partial [\sigma^{-1}]_{ij}(x)}{\partial x_k} = \frac{\partial}{\partial x_k}\left(\frac{\partial \gamma_i(x)}{\partial x_j}\right) = \frac{\partial}{\partial x_j}\left(\frac{\partial \gamma_i(x)}{\partial x_k}\right) = \frac{\partial [\sigma^{-1}]_{ik}(x)}{\partial x_j},$$

for all $(i,j,k) = 1,\ldots,m$. Continuity of the second-order partial derivatives is required for the order of differentiation to be interchangeable. Here, we have infinite differentiability.

Conversely, suppose that $\sigma^{-1}$ satisfies (10). Then, for each $i = 1,\ldots,m$, use row $i$ of the matrix $\sigma^{-1}$, $\sigma_{i\cdot}^{-1} = [[\sigma^{-1}]_{ij}]_{j=1,\ldots,m}$, to define the differential 1-form $\omega_i = \sum_{j=1}^{m} [\sigma^{-1}]_{ij}\, dx_j$ and calculate its differential, the differential 2-form $d\omega_i$. Condition (10) implies that $d\omega_i = 0$, that is the differential 1-form $\omega_i$ is closed on $\mathcal{S}_X$. The domain $\mathcal{S}_X$ is singly connected (or without holes). Therefore, by Poincaré's lemma (see, e.g., Theorem V.8.1 of Edwards [12]), the form $\omega_i$ is exact, that is there exists a differential 0-form $\gamma_i$ such that $d\gamma_i = \omega_i$. In other words, for each row $i$ of the matrix $\sigma^{-1}$, there exists a function $\gamma_i$ defined by $\gamma_i(x) = \int^{x_j} [\sigma^{-1}]_{ij}(x)\, dx_j$ (the choice of the index $j$ is irrelevant) which satisfies (50), has the required differentiability properties and is invertible. The function $\gamma$ is then defined by each of its $d$ components $\gamma_i$, $i = 1,\ldots,m$, and because of Assumptions 2 and 3, it is invertible and infinitely often differentiable. By construction, $Y_t \equiv \gamma(X_t)$ has unit diffusion and therefore $X$ is reducible. To prove the equivalent characterization (9), apply Itô's lemma from $Y$ to $X$ (instead of from $X$ to $Y$) and proceed as above.

8.2. *Proof of Theorem* 1. To show that (15) with the coefficients given in the statement of Theorem 1 indeed represent the Taylor expansion in $\Delta$ of the log-density function $l_Y$, at order $K-1$, it suffices to verify that the difference between the left- and right-hand sides in the Kolmogorov forward and backward partial differential equations is of order $\Delta^K$.

Define $F_Y^{(K)}(y|y_0,\Delta)$ [resp. $B_Y^{(K)}(y|y_0,\Delta)$] as the difference between the left- and right-hand sides of the forward (resp. backward) equations when $l_Y$ is replaced by $l_Y^{(K)}$. The backward equation for $l_Y$ is

$$
\begin{aligned}
\frac{\partial l_Y(y|y_0,\Delta)}{\partial \Delta} &= \sum_{i=1}^{m} \mu_{Yi}(y_0) \frac{\partial l_Y(y|y_0,\Delta)}{\partial y_{0i}} \\
&\quad + \frac{1}{2}\sum_{i=1}^{m} \frac{\partial^2 l_Y(y|y_0,\Delta)}{\partial y_{0i}^2} + \frac{1}{2}\sum_{i=1}^{m}\left(\frac{\partial l_Y(y|y_0,\Delta)}{\partial y_{0i}}\right)^2.
\end{aligned}
\tag{51}
$$

Substituting in the expansion (15) and equating the coefficients of $\Delta^{-2}$ on both sides of (51) yields

$$C_Y^{(-1)}(y|y_0) = -\frac{1}{2}\left(\frac{\partial C_Y^{(-1)}(y|y_0)}{\partial y_{0i}}\right)^T\left(\frac{\partial C_Y^{(-1)}(y|y_0)}{\partial y_{0i}}\right),$$

which is satisfied by the (already determined) solution (20), which is therefore the desired solution.

Starting with the Gaussian leading term (20), we have

$$\begin{cases} F_Y^{(K)}(y|y_0,\Delta) = \displaystyle\sum_{k=-1}^{K-1} f_Y^{(k)}(y|y_0)\frac{\Delta^k}{k!} + O(\Delta^K) \\ B_Y^{(K)}(y|y_0,\Delta) = \displaystyle\sum_{k=-1}^{K-1} b_Y^{(k)}(y|y_0)\frac{\Delta^k}{k!} + O(\Delta^K) \end{cases}$$

[with the convention that $(-1)! = 0! = 1$]. The first term in $F_Y^{(K)}$ is

$$\begin{aligned}f_Y^{(-1)}(y|y_0) = {}& -\frac{m}{2} + \sum_{i=1}^m \mu_{Yi}(y)\frac{\partial C_Y^{(-1)}(y|y_0)}{\partial y_i} - \frac{1}{2}\sum_{i=1}^m \frac{\partial^2 C_Y^{(-1)}(y|y_0)}{\partial y_i^2} \\ & -\frac{1}{2}\sum_{i=1}^m 2\frac{\partial C_Y^{(-1)}(y|y_0)}{\partial y_i}\frac{\partial C_Y^{(0)}(y|y_0)}{\partial y_i} \\ ={}& -\sum_{i=1}^m (y_i - y_{0i})\mu_{Yi}(y) + \sum_{i=1}^m (y_i - y_{0i})\frac{\partial C_Y^{(0)}(y|y_0)}{\partial y_i}.\end{aligned}$$

Solving the equation $f_Y^{(-1)}(y|y_0) = 0$ for $C_Y^{(0)}(y|y_0)$ yields the full solution

(52)
$$\begin{aligned} C_Y^{(0)}(y|y_0) = {}& \sum_{i=1}^m (y_i - y_{0i})\int_0^1 \mu_{Yi}(y_0 + u(y-y_0))\,du \\ & + \sum_{i,j=1,\ j\neq i}^m \alpha_{ij}^{(0)}\frac{y_i - y_{0i}}{y_j - y_{0j}} + M_Y^{(0)}, \end{aligned}$$

where the $\alpha_{ij}^{(0)}$ and $M_Y^{(0)}$ are integration constants in the differential equation $f_Y^{(-1)} = 0$, hence arbitrary functions of $y_0$. The boundary condition that $C_Y^{(0)}$ be finite when passing through the axes $y_j = y_{j0}$ for all $j = 1,\ldots,m$ imposes the condition $\alpha_{ij}^{(0)} = 0$. To determine $M_Y^{(K)}$, (51) gives

$$b_Y^{(-1)}(y|y_0) = -\sum_{i=1}^m (y_i - y_{0i})\mu_{Yi}(y) - \sum_{i=1}^m (y_i - y_{0i})\frac{\partial C_Y^{(0)}(y|y_0)}{\partial y_{0i}}$$

and the candidate solution (52) must satisfy $b_Y^{(-1)} = 0$. Thus, we must have

$$\sum_{i=1}^{m} \frac{\partial M_Y^{(0)}(y_0)}{\partial y_{0i}} (y_i - y_{0i}) = 0$$

for all $y$ and $y_0$ and so $M_Y^{(0)}(y_0)$ is constant. Since the limiting behavior of $p_Y$ is $N(0, I)$ as $\Delta \to 0$ and

$$\lim_{\Delta \to 0} \left( l_Y^{(K)}(y|y_0, \Delta) + \frac{m}{2} \ln(2\pi\Delta) + \frac{1}{2\Delta} \|y - y_0\|^2 \right) = C_Y^{(0)}(y|y_0),$$

we must have $M_Y^{(0)} = 0$ to ensure that as $\Delta \to 0$, the limiting density integrates to one [otherwise, it integrates to $\exp(M_Y^{(0)})$].

The next term is

$$f_Y^{(0)}(y|y_0) = C_Y^{(1)}(y|y_0) + \sum_{i=1}^{m}(y_i - y_{0i})\frac{\partial C_Y^{(1)}(y|y_0)}{\partial y_i} + \sum_{i=1}^{m} \frac{\partial \mu_{Yi}(y)}{\partial y_i}$$

$$+ \sum_{i=1}^{m} \mu_{Yi}(y) \frac{\partial C_Y^{(0)}(y|y_0)}{\partial y_i}$$

$$- \frac{1}{2} \sum_{i=1}^{m} \left\{ \frac{\partial^2 C_Y^{(0)}(y|y_0)}{\partial y_i^2} + \left[ \frac{\partial C_Y^{(0)}(y|y_0)}{\partial y_i} \right]^2 \right\}$$

$$= C_Y^{(1)}(y|y_0) + \sum_{i=1}^{m}(y_i - y_{0i})\frac{\partial C_Y^{(1)}(y|y_0)}{\partial y_i} - G_Y^{(1)}(y|y_0),$$

where $G_Y^{(1)}$ is given in (23) and depends on the previously determined $C_Y^{(-1)}$ and $C_Y^{(0)}$.

Solving the equation $f_Y^{(0)}(y|y_0) = 0$, which is linear in $C_Y^{(1)}$, similarly yields the explicit solution

$$C_Y^{(1)}(y|y_0) = \int_0^1 G_Y^{(1)}(y_0 + u(y - y_0)|y_0)\,du$$

$$+ \sum_{i,j=1, j\neq i}^{m} \alpha_{ij}^{(1)} \frac{y_i - y_{0i}}{(y_j - y_{0j})^2} + M_Y^{(1)},$$

which includes generic integration constants $\alpha_{ij}^{(1)}$ and $M_Y^{(1)}$. The solution has $\alpha_{ij}^{(1)} = 0$ when imposing finiteness of $l_Y^{(K)}$ when passing through the axes $y_j = y_{j0}$ for all $j = 1, \ldots, m$. As for $M_Y^{(1)}$, invoking the backward equation (51) yields

$$M_Y^{(1)}(y_0) - \sum_{i=1}^{m} \frac{\partial M_Y^{(1)}(y_0)}{\partial y_{0i}}(y_i - y_{0i}) = 0,$$

whose only solution valid for all $(y, y_0)$ is $M_Y^{(1)}(y_0) = 0$. This yields the coefficient $C_Y^{(1)}$.

More generally, the term $f_Y^{(k-1)}$, $k \geq 1$, is given by

$$f_Y^{(k-1)}(y|y_0) = C_Y^{(k)}(y|y_0) + \frac{1}{k} \sum_{i=1}^{m} (y_i - y_{0i}) \frac{\partial C_Y^{(k)}(y|y_0)}{\partial y_i} - G_Y^{(k)}(y|y_0),$$

where $G_Y^{(k)}$ is given in (24) and depends on the previously determined $C_Y^{(-1)}, C_Y^{(0)}, \ldots, C_Y^{(k-1)}$. Solving the equation $f_Y^{(k)}(y|y_0) = 0$ for $C_Y^{(k)}$ (with the same boundary condition as for $C_Y^{(0)}$ and $C_Y^{(1)}$) yields the explicit solution (22). In this case, the full solution including generic integration constants $\alpha_{ij}$ and $M_Y^{(k)}$, is

$$C_Y^{(k)}(y|y_0) = k \int_0^1 G_Y^{(k)}(y_0 + u(y - y_0)|y_0) u^{k-1}\, du$$

$$+ \sum_{i,j=1, j \neq i}^{m} \alpha_{ij}^{(k)} \frac{y_i - y_{0i}}{(y_j - y_{0j})^{k+1}} + M_Y^{(k)}.$$

Thus, by construction, the solution $C_Y^{(k)}$, $k = -1, 0, \ldots, K$ given in the statement of the theorem is such that $F_Y^{(K)}(y|y_0, \Delta) = O(\Delta^K)$. Similarly, $B_Y^{(K)}(y|y_0, \Delta) = O(\Delta^K)$. The fact that solving the Kolmogorov equations to order $\Delta^K$ yields a Taylor expansion of order $K-1$ of $l_Y$ is established as part of the proof of Theorem 3 below.

8.3. *Proof of Theorem* 2. Let $F_X^{(K)}$ and $\tilde{F}_X^{(K)}$ (resp. $B_X^{(K)}$ and $\tilde{B}_X^{(K)}$) denote the difference between the left- and right-hand sides of (28) [resp. (29)] when $l_X$ is replaced by the expansion $l_X^{(K)}$ (resp. $\tilde{l}_X^{(K)}$). We have

$$F_X^{(K)}(x|x_0, \Delta) = \sum_{k=-2}^{K-1} f_X^{(k)}(x|x_0) \frac{\Delta^k}{k!} + O(\Delta^K)$$

[with the convention that $(-2)! = 2$ and $(-1)! = 0! = 1$]. The highest-order term is $f_X^{(-2)}$, given by (39), and the coefficient function $C_X^{(-1)}$ is such that it sets $f_X^{(-2)}$ to zero. We have then successively determined $C_X^{(0)}$ by setting $f_X^{(-1)}$ in (40) to zero and, more generally, given $C_X^{(-1)}, C_X^{(0)}, \ldots, C_X^{(k-1)}$, the expression (41) for $f_X^{(k-1)}$ is defined and can be set to zero to determine the next coefficient $C_X^{(k)}$. The form of the log-likelihood adopted in (27) with $D_v$ kept separate from $C_X^{(0)}$ is essential to obtain $\tilde{B}_X^{(K)}(x|x_0, \Delta) = O(\Delta^K)$ in addition to $\tilde{F}_X^{(K)}(x|x_0, \Delta) = O(\Delta^K)$.

To determine the expansions in $x - x_0$ for each coefficient $C_X^{(k)}$, $k \geq -1$, replace $C_X^{(k)}$ by $C_X^{(j_k, k)}$ in each equation in turn. Starting with (39), calculate an expansion of $\tilde{f}_X^{(-2)}$ in $(x - x_0)$ to order $j_{-1}$. This determines a system of equations in the unknown coefficients $\beta_i^{(-1)}$, $i \in I_{-1}$ [which appear when $C_X^{(-1)}$ is expanded, as in (34)]. By construction, there are as many equations as unknowns (both are given by the number of elements in $I_{-1}$). This system of equations can always be solved explicitly because it has the following form. First, $\beta_i^{(-1)} = 0$ for $\mathrm{tr}[i] = 0, 1$ (i.e., the polynomial has no constant or linear terms) and the terms corresponding to $\mathrm{tr}[i] = 2$ (with, of course, $j_{-1} \geq 2$) are

$$\sum_{i \in I_{-1}: \mathrm{tr}[i]=2} \beta_i^{(-1)}(x_0)(x_1 - x_{01})^{i_1} \cdots (x_m - x_{0m})^{i_m}$$
$$= -(x - x_0)^T v^{-1}(x_0)(x - x_0),$$

which is the anticipated term, given the Gaussian limiting behavior of the transition density when $\Delta$ is small. Thus, with $j_{-1} \geq 3$, we only need to determine the terms $\beta_i^{(-1)}$ corresponding to $\mathrm{tr}[i] = 3, \ldots, j_{-1}$. Then, the next order coefficients in $(x - x_0)$, that is, the coefficients corresponding to $\mathrm{tr}[i] = 3$, each appear linearly in a separate equation. That is, we have a system $\Phi_3^{(-1)}(x_0) \cdot \beta_3^{(-1)}(x_0) = a_3^{(-1)}(x_0)$ whose explicit solution is given by $\beta_3^{(-1)}(x_0) = \mathrm{Inv}[\Phi_3^{(-1)}(x_0)] \cdot a_3^{(-1)}(x_0)$. Given the previously determined coefficients corresponding to $\mathrm{tr}[i] = 0, \ldots, r$, the equations determining the coefficients for $\mathrm{tr}[i] = r + 1$ are given by a linear system $\Phi_{r+1}^{(-1)}(x_0) \cdot \beta_{r+1}^{(-1)}(x_0) = a_{r+1}^{(-1)}(x_0)$, where the matrix $\Phi_{r+1}^{(-1)}$ and the vector $a_{r+1}^{(-1)}$ are functions of the previously determined coefficients $\beta_i^{(-1)}$, for $\mathrm{tr}[i] = 0, \ldots, r$, and $x_0$.

The same principle applies to all values of $k$. For $k = 0$, $\beta_i^{(0)} = 0$ for $\mathrm{tr}[i] = 0$, so the polynomial has no constant term. For $k \geq 1$, the polynomials have a constant term (for $k \geq 1$, $\beta_i^{(k)} \neq 0$ for $\mathrm{tr}[i] = 0$, in general). The same principle applies to each equation in turn: once $C_X^{(j_{-1}, -1)}$ is determined, a Taylor expansion of (40) determines the coefficients $\beta_i^{(0)}$, $i \in I_0$, and so on.

8.4. *Proof of Proposition* 2. If the diffusion $X$ is reducible, then $C_X^{(k)}(x|x_0) = C_Y^{(k)}(\gamma(x)|\gamma(x_0))$. By construction (see the proof of Theorem 2), the coefficients $C_X^{(j_k, k)}$ are then Taylor expansions of the coefficients $C_X^{(k)}$ (which are the solutions of the equations $f_X^{(k-1)} = 0$).

8.5. *Proof of Theorem* 3. Consider the irreducible case; everything also applies to the reducible case, by simply eliminating the arguments relative to the additional expansion in $x - x_0$. The expansion $\tilde{p}_X^{(K)}$, constructed as an expansion in $\Delta$ and $(x - x_0)$ of $\exp(\tilde{l}_X^{(K)})$, satisfies the linear backward equation, but with a remainder term $\tilde{B}_X^{(K)}$:

$$
\begin{aligned}
\tilde{B}_X^{(K)}(x|x_0, \Delta) &= \frac{\partial \tilde{p}_X^{(K)}(x|x_0, \Delta)}{\partial \Delta} - \sum_{i=1}^m \mu_i(x_0) \frac{\partial \tilde{p}_X^{(K)}(x|x_0, \Delta)}{\partial y_{0i}} \\
&\quad - \frac{1}{2} \sum_{i,j=1}^m v_{ij}(x_0) \frac{\partial^2 \tilde{p}_X^{(K)}(x|x_0, \Delta)}{\partial x_{0i} \partial x_{0j}}.
\end{aligned}
\tag{53}
$$

Indeed, the coefficients $C_X^{(j_k, k)}$, $k = -1, \ldots, K$, are constructed in such a way that the terms of order $\Delta^{K-1}$ and higher of the right-hand side of (53) are zero. The remainder term is of the form $\tilde{B}_X^{(K)}(x|x_0, \Delta) = \Delta^K \tilde{p}_X^{(K)}(x|x_0, \Delta) \tilde{\psi}_X^{(K)}(x|x_0, \Delta)$ where $\tilde{\psi}_X^{(K)}(x|x_0, \Delta)$ is a sum of products of the functions $\mu_i$, $v_{ij}$, the coefficients $C_X^{(j_k, k)}$, $k = -1, \ldots, K$, and their derivatives. Specifically, $\tilde{\psi}_X^{(K)}(x|x_0, \Delta) = c_K f_X^{(K)}(x|x_0) + o(\Delta)$, where the functions $f_X^{(k)}$ are given in Theorem 2 and where $c_K$ is a constant.

The coefficients and their derivatives exhibit at most polynomial growth as a result of the explicit expressions given in Sections 4 and 5, combined with Assumption 4 on $(\mu, v)$ and their derivatives. Thus, $\tilde{\psi}_X^{(K)}$ exhibits at most polynomial growth and we have $\tilde{\psi}_X^{(K)}(x|x_0, \Delta) = O(1)$ uniformly for all $(x, x_0)$ in a compact subset of the interior of $\mathcal{S}_X^2$ and for all $\theta$ in $\Theta$, by virtue of the continuity of the functions and their derivatives with respect to the parameter vector in the compact set $\Theta$.

The solution $\tilde{p}_X^{(K)}$ of the approximate PDE with remainder, (53), is an approximation of the solution $p_X$ of the exact PDE without remainder term due to the following. Writing $\tilde{r}_X^{(K)} = \tilde{p}_X^{(K)} - p_X$, it is clear, by linearity of the PDE for $p_X$, that $\tilde{r}_X^{(K)}$ also satisfies equation (53) with the same remainder $\tilde{B}_X^{(K)}$, but now with initial condition $\tilde{r}_X^{(K)}(x|x_0, \Delta) \to 0$ as $\Delta \to 0$ (whereas $p_X$ and $\tilde{p}_X^{(K)}$ both converge to a Dirac mass at $x_0$ as $\Delta \to 0$). The solution is given by

$$
\tilde{r}_X^{(K)}(x|x_0, \Delta) = \int_{\mathcal{S}_X} \int_0^\Delta \tilde{B}_X^{(K)}(x|z, \tau) p_X(z|x_0, \Delta - \tau) \, d\tau \, dz,
\tag{54}
$$

which follows from the facts that (54) produces the correct initial boundary behavior and, as can be seen by computing the relevant partial derivatives of this expression for $\tilde{r}_X^{(K)}$, satisfies (53). The function $\tilde{B}_X^{(K)} p_X$ is integrable because $p_X$ has exponentially decaying tails in a neighborhood of $\Delta = 0$ (see below), whereas $\phi_X^{(K)}$ has polynomial growth. It follows that we have $\tilde{r}_X^{(K)}(x|x_0, \Delta) = O(\Delta^K)$ uniformly for all $(x, x_0)$ in a compact subset of the interior of $\mathcal{S}^2$ and for all $\theta$ in $\Theta$. Let

$$\tilde{R}_X^{(K)}(x|x_0, \Delta) \equiv \sup_{\theta \in \Theta} |\tilde{r}_X^{(K)}(x|x_0, \Delta)|$$

and consider the expectation

$$(55) \quad E[\tilde{R}_X^{(K)}(X_\Delta|X_0, \Delta)|X_0 = x_0] = \int_{\mathcal{S}_X} \tilde{R}_X^{(K)}(x|x_0, \Delta) p_X(x|x_0, \Delta) \, dx.$$

In the integral above, divide the region $\mathcal{S}_X \subseteq \mathbb{R}^m$ into two parts—a neighborhood of $x_0$ of the form $\mathcal{N} = \prod_{i=1}^m [x_{0i} - \Delta^{1/2} c_\Delta, x_{0i} + \Delta^{1/2} c_\Delta]$, where $c_\Delta$ is a sequence of positive numbers such that $c_\Delta \to \infty$ and $\Delta^{1/2} c_\Delta \to 0$, and the complement $\mathcal{S}_X \setminus \mathcal{N}$ of that neighborhood.

From above, there exist constants $C$ and $M$ such that $|\tilde{R}_X^{(K)}(x|x_0, \Delta)| \leq M \Delta^K$ for all $x \in \mathcal{N}$, so

$$(56) \quad \int_\mathcal{N} \tilde{R}_X^{(K)}(x|x_0, \Delta) p_X(x|x_0, \Delta) \, dx = O(\Delta^K).$$

Outside the neighborhood $\mathcal{N}$, the approximation error $\tilde{R}_X^{(K)}$ need not be small; however, it is at most polynomial in $(x - x_0)$. In a neighborhood of $\Delta = 0$, the tail behavior of the transition function $p_X$ is driven by the term $\exp(-\frac{m}{2} \ln(\Delta) + C_X^{(j-1,-1)}(x|x_0) \Delta^{-1})$, with $C_X^{(j-1,-1)}(x|x_0) = -(1/2)(x - x_0)^T v^{-1}(x_0)(x - x_0) + o(\|x - x_0\|^2)$.

Therefore, to bound $\int_{\mathcal{S}_X \setminus \mathcal{N}} \tilde{r}_X^{(K)} p_X \, dx$, one needs to integrate a polynomial error term in $r_X^{(K)}$, say $\|x - x_0\|^b$, $b \geq 0$, against the exponential tails of $p_X$. This results in integrals of the form (written in the univariate case $m = 1$, for simplicity, and near an infinity boundary)

$$(57) \quad \begin{aligned} &\Delta^{-1/2} \int_{\Delta^{1/2} c_\Delta}^{+\infty} |x - x_0|^b \exp(-(x - x_0)^2/(2\Delta v(x_0))) \, dx \\ &= \Delta^{b/2} \int_{c_\Delta}^{+\infty} |z - z_0|^b \exp(-(z - z_0)^2/(2 v(x_0))) \, dz, \end{aligned}$$

with the change of variable $z - z_0 = (x - x_0)/\Delta^{1/2}$, and similarly on the interval $[-\infty, -\Delta^{1/2} c_\Delta]$. Since $c_\Delta \to \infty$, the above integral converges to zero.

It follows from (56) and (57) that $E[\tilde{R}_X^{(K)}(X_\Delta|X_0, \Delta)|X_0] \to 0$ as $\Delta \to 0$. Similar calculations show that this is also the case for $\mathrm{Var}[\tilde{R}_X^{(K)}(X_\Delta|X_0, \Delta)|X_0]$. Therefore, $\tilde{R}_X^{(K)}(X_\Delta|X_0, \Delta) \to 0$ in $P_{\theta_0}$-probability, given $X_0$, which follows by Chebyshev's inequality from the convergence to zero of the conditional expected value and variance of $\tilde{R}_X^{(K)}$. Convergence to zero of the conditional probability implies convergence to zero of the unconditional probability, since

$$\Pr(|\tilde{R}_X^{(K)}(X_{t+\Delta}|X_t, \Delta)| > \varepsilon)$$
$$= \int_{\mathcal{S}_X} \Pr(|\tilde{R}_X^{(K)}(X_{t+\Delta}|X_t, \Delta)| > \varepsilon | X_t = x_0) \pi_t(x_0) \, dx_0,$$

where $\pi_t(x_0)$ denotes the marginal density of $X$ at time $t$. Since probabilities are between 0 and 1 and since $\pi_t$ integrates to 1, the convergence of the unconditional probability follows by Lebesgue's dominated convergence theorem. Next, the convergence in $P_{\theta_0}$-probability for the log-density follows from the continuity of the logarithm and the convergence of $\tilde{R}_X^{(K)}$. Then, for fixed $n$, the convergence stated in (42) follows from that of $\tilde{l}_X^{(K)}$ to $l_X$. From the assumed existence of the maximizer $\hat{\theta}_{n,\Delta}$ of $\tilde{\ell}_n(\theta, \Delta)$ in $\Theta$ and the (just obtained) proximity of the two objective functions, it follows from standard arguments that the maximizer of $\tilde{\ell}_n^{(K)}(\theta, \Delta)$, that is, $\hat{\theta}_{n,\Delta}^{(K)}$, exists almost surely, is in $\Theta$ and is close to $\hat{\theta}_{n,\Delta}$ as $\Delta \to 0$, in the sense that $\hat{\theta}_{n,\Delta}^{(K)} - \hat{\theta}_{n,\Delta} \to 0$ in $P_{\theta_0}$-probability. Finally, the speed at which $\hat{\theta}_{n,\Delta}^{(K)} - \hat{\theta}_{n,\Delta}$ converges to zero can be made arbitrarily small for any $n$ by taking $\Delta \to 0$ sufficiently fast. In particular, as $\to \infty$, a sequence $\Delta_n \to 0$ can be taken to be such that (44) is satisfied.

**9. Conclusions.** This paper provides a method to derive closed-form Taylor expansions to the likelihood function of arbitrary multivariate diffusions. While these expansions are local in nature, at least in the irreducible case, they have been shown to produce useful approximations in the context of maximum likelihood estimation. Aït-Sahalia and Kimmel [3] apply this method to popular stochastic volatility models. Likelihood expansions for these models are derived, as well as a Monte Carlo investigation of the properties of maximum likelihood estimators of the parameters computed from these expansions.

Once the expansion is computed for the diffusion model at hand, it can be applied to the estimation of parameters by a variety of other estimation methods which require an expression for the transition density of the state variables, such as Bayesian methods where one wishes to obtain a posterior distribution for the parameters of a stochastic differential equation, or to generate simulated data at the desired frequency from the continuous-time model or to serve as the instrumental or auxiliary model in indirect inference and simulated or efficient moments methods. The explicit nature of the expansion as a function of all of the relevant variables makes these computations, whether maximization of the classical likelihood or computation of posterior distributions, straightforward and computationally very efficient.

Other methods can be used to approximate the transition function: they involve numerically solving a partial differential equation, simulating the process to Monte Carlo integrate the transition density or approximating the process with binomial trees. However, none of these alternative methods provides closed-form formulae. Jensen and Poulsen [20], Stramer and Yan [26] and Hurn, Jeisman and Lindsay [19] compare the accuracy and speed of the different methods, showing the advantages of this closed-form approach.

**Acknowledgments.** I am very grateful to the Editor, the Associate Editor and three referees for constructive and helpful comments and suggestions.

DEPARTMENT OF ECONOMICS
PRINCETON UNIVERSITY
AND
NBER
PRINCETON, NEW JERSEY 08544-1021
USA
E-MAIL: yacine@princeton.edu